\newcommand{\al}{\alpha}

\newcommand{\de}{\delta}               

\newcommand{\lb}{\lambda}

\newcommand{\vphi}{\varphi}

\newcommand{\cal}{\mathcal}

\newcommand{\calf}{{\cal F}}

\newcommand{\calk}{{\cal K}}

\newcommand{\calr}{{\cal R}}           

\newcommand{\calt}{{\cal T}}           
          
\newcommand{\calv}{{\cal V}}

\newcommand{\Fix}{{\rm Fix}}

\newcommand{\incl}{\subseteq}

\newcommand{\es}{\emptyset}          
\newcommand{\sm}{\setminus}
      
\newcommand{\limpl}{\Longrightarrow}

\newcommand{\oo}{\infty}

\newcommand{\sk}{\smallskip}

\def\dtends   {\stackrel {\it d}{\longrightarrow}}

\def\Dtends   {\stackrel {\it D}{\longrightarrow}}

\def\(V)tends {\stackrel {(\calv)}{\longrightarrow}}

\newcommand{\barr}{\begin{array}}        
\newcommand{\earr}{\end{array}}
\newcommand{\bcor}{\begin{corollary}}    
\newcommand{\ecor}{\end{corollary}}
\newcommand{\ben}{\begin{enumerate}}     
\newcommand{\een}{\end{enumerate}}
\newcommand{\beq}{\begin{equation}}       
\newcommand{\eeq}{\end{equation}}
\newcommand{\bex}{\begin{example}}        
\newcommand{\eex}{\end{example}}
\newcommand{\bit}{\begin{itemize}}        
\newcommand{\eit}{\end{itemize}}
\newcommand{\blemma}{\begin{lemma}}       
\newcommand{\elemma}{\end{lemma}}
\newcommand{\bproof}{\begin{proof}}       
\newcommand{\eproof}{\end{proof}}
\newcommand{\bprop}{\begin{proposition}}  
\newcommand{\eprop}{\end{proposition}}
\newcommand{\brem}{\begin{remark}}        
\newcommand{\erem}{\end{remark}}
\newcommand{\btab}{\begin{tabular}}       
\newcommand{\etab}{\end{tabular}}
\newcommand{\btheorem}{\begin{theorem}}   
\newcommand{\etheorem}{\end{theorem}}

\documentclass[reqno]{amsart}

\newtheorem{theorem}{\bf Theorem}
\newtheorem{corollary}{\bf Corollary}
\newtheorem{example}{\bf Example}
\newtheorem{lemma}{\bf Lemma}
\newtheorem{proposition}{\bf Proposition}
\newtheorem{remark}{\bf Remark}

\begin{document}

\title
[PPF Dependent Fixed Points in A-closed Razumikhin Classes]
{ PPF DEPENDENT FIXED POINTS IN \\
A-CLOSED RAZUMIKHIN CLASSES}

\author{Mihai Turinici}
\address{
"A. Myller" Mathematical Seminar;
"A. I. Cuza" University;
700506 Ia\c{s}i, Romania
}
\email{mturi@uaic.ro}




\begin{abstract}
The PPF dependent fixed point result in
algebraically closed Razumikhin classes due to
Agarwal et al 
[Fixed Point Theory Appl., 2013, 2013:280]
is identical with its 
constant class counterpart;
and this, in turn, is
reducible to a fixed point principle
involving SVV type contractions
(over the subsequent metric space
of initial Banach structure),
without any regularity conditions 
about the Razumikhin classes. 
The conclusion remains valid 
for all PPF dependent fixed point results 
founded on such global conditions.  
\end{abstract}

\maketitle

\section{Introduction}
\setcounter{equation}{0}

Let $X$ be a nonempty set.
Call the subset $Y$ of $X$, 
{\it almost singleton} (in short: {\it asingleton})
provided $y_1,y_2\in Y$ implies $y_1=y_2$;
and {\it singleton},
if, in addition, $Y$ is nonempty;
note that, in this case,
$Y=\{y\}$, for some $y\in X$. 
Take a metric $d:X\times X\to R_+:=[0,\oo[$ over it;
as well as a selfmap 
$T\in \calf(X)$.
[Here, for each couple $A,B$ of nonempty sets,
$\calf(A,B)$ denotes the class of all functions 
from $A$ to $B$; 
when $A=B$, we write $\calf(A)$ in place of $\calf(A,A)$].
Denote $\Fix(T)=\{x\in X; x=Tx\}$;
each point of this set is referred to as 
{\it fixed} under $T$.
The determination of such elements is to be performed
under the directions below,
comparable with the ones in 
Rus \cite[Ch 2, Sect 2.2]{rus-2001}
and
Turinici \cite{turinici-2011-JIMS}:

{\bf Pic-1)} 
We say that $T$ is a {\it Picard operator} (modulo $d$) if,
for each $x\in X$, the iterative sequence
$(T^nx; n\ge 0)$ is $d$-convergent; 
and a {\it globally Picard operator} (modulo $d$) 
if, in addition, 
$\Fix(T)$ is an asingleton

{\bf Pic-2)} 
We say that $T$ is a {\it strong Picard operator} (modulo $d$) if,
for each $x\in X$, the iterative sequence
$(T^nx; n\ge 0)$ is $d$-convergent 
with $\lim_n(T^nx)$ belonging to $\Fix(T)$;
and a {\it globally strong Picard operator} (modulo $d$) 
if, in addition, 
$\Fix(T)$ is an asingleton (hence, a singleton).

The basic result in this area is the 1922 one due to
Banach \cite{banach-1922}.
Given $k\ge 0$, let us say that  $T$ is  
{\it $(d;k)$-contractive}, if;
\ben
\item[] (a01)\ \ 
$d(Tx,Ty)\le k d(x,y)$,\ for all  $x,y\in X$.
\een

\btheorem \label{t1}
Assume that $T$ is $(d;k)$-contractive, for some $k\in [0,1[$.
In addition, let $(X,d)$ be complete.
Then, 
$T$ is globally strong Picard (modulo $d$). 
\etheorem

This statement 
(referred to as: Banach's fixed point principle) 
found some basic applications to 
different branches of
operator equations theory.
As a consequence, many 
extensions of it were proposed.
From the perspective of this exposition,
the following ones are of interest:

{\bf I)}
Contractive type extensions:
the initial Banach contractive property is 
taken in a generalized way, as
\ben
\item[] (a02)\ \ 
$F(d(Tx,Ty),d(x,y),d(x,Tx),d(y,Ty),d(x,Ty),d(Tx,y))\le 0$,\\
for all $x,y\in X$;
\een
where $F:R_+^6\to R$ is a function.
For the explicit case, 
some consistent lists of these 
may be found in 
the survey papers by 
Rhoades \cite{rhoades-1977},
Collaco and E Silva \cite{collaco-e-silva-1997},
Kincses and Totik \cite{kincses-totik-1999},
as well as the references therein.
And, for the implicit case, 
certain particular aspects 
have been considered by 
Leader \cite{leader-1979}
and
Turinici \cite{turinici-1976-AUAICI}.

{\bf II)}
Structural extensions:
the trivial relation $X\times X$ 
is replaced by a 
relation $\nabla$ over $X$ fulfilling or not
certain regularity properties.
For example, the case of 
$\nabla$ being  
{\it reflexive, transitive}
(hence, a {\it quasi-order})
was considered in the 1986 paper by
Turinici \cite{turinici-1986-DM}.
Two decades later,
this result was re-discovered --
in a (partially) ordered context --
by
Ran and Reurings \cite{ran-reurings-2004};
see also 
Nieto and Nodriguez-Lopez \cite{nieto-rodriguez-lopez-2005}.
On the other hand, the 
"amorphous" case 
($\nabla$ has no regularity properties at all)
was discussed (via graph techniques) in 
Jachymski \cite{jachymski-2008};
and (from a general perspective) by
Samet and Turinici \cite{samet-turinici-2012}.
Some other aspects involving 
additional convergence structures 
may be found in 
Kasahara \cite{kasahara-1976}.

{\bf III)}
Nonself extensions:
$T$ is no longer a selfmap.
In 1997, 
Bernfeld et al 
\cite{bernfeld-lakshmikantham-reddy-1977} 
introduced the concept of 
{\it PPF (past-present-future) dependent fixed point}
for nonself-mappings (whose domain is distinct from their range).
Furthermore, the quoted authors established --
via iterative methods involving
a certain Razumikhin class $\calr_c$ --
some PPF dependent fixed point theorems 
for contractive mappings of this type. 
As precise there, 
the obtained results are  
useful tools in the study of existence and uniqueness questions 
for solutions of nonlinear functional differential/integral equations 
which may depend upon
the past history, present data and future evolution.
As a consequence, this theory attracted a lot of
contributors in the area;
see, for instance,
Dhage \cite{dhage-2012-FPT, dhage-2012-JNSA},
Hussain et al \cite{hussain-khaleghizadeh-salimi-akbar-2013},
Kaewcharoen \cite{kaewcharoen-2013},
Kutbi and Sintunavarat \cite{kutbi-sintunavarat-2014},
as well as the references therein.
However, as proved in a recent paper by
Cho, Rassias, Salimi and Turinici 
\cite{cho-rassias-salimi-turinici-2014},
the starting conditions
[imposed by the problem setting]
relative to the 
ambient Razumikhin class $\calr_c$
may be converted into  
starting conditions relative to 
the constant class $\calk$;
so, ultimately, 
we may arrange for these PPF dependent fixed point
results holding over $\calk$. 
In this exposition we 
bring the discussion a step
further, by establishing that

{\bf Fact-1)}  
the algebraic closeness assumption 
[used in all these references]
imposed upon the Razumikhin class $\calr_c$
yields, in a direct way, $\calr_c=\calk$

{\bf Fact-2)}
the PPF dependent fixed point problem
attached to the constant class $\calk$
is reducible to a (standard) fixed point problem 
in the (complete) metrical structure
induced by our initial Banach one,
under no regularity assumption about
the ambient Razumikhin class.

Finally, as an application 
of these conclusions,
we show that a recent PPF dependent 
fixed point result in
Agarwal et al \cite{agarwal-kumam-sintunavarat-2013}
is reducible to a fixed point problem 
involving a class of contractions
over standard metric structures 
[taken as before]
introduced under the lines proposed by
Samet et al \cite{samet-vetro-vetro-2012}.
Further aspects will be discussed elsewhere.

\section{Razumikhin classes}
\setcounter{equation}{0}

Let $(E,||.||)$ be a Banach space;
and $d(.,.)$ be the induced by norm metric on $E$:
\ben
\item[]
$d(x,y)=||x-y||$,\ $x,y\in E$;
\een
hence, $(E,d)$ is a complete metric space.
Further, let $I=[a,b]$ be a closed real interval,
and $E_0:=C(I,E)$ stand for the class 
of all continuous functions 
$\vphi:I\to E$, endowed with the supremum norm
\ben
\item[] 
$||\vphi||_{0}=\sup\{||\vphi(t)||; t\in I\}$,\ $\vphi\in E_0$.
\een
As before, let $D(.,.)$ stand for the 
induced by norm metric on $E_0$
\ben
\item[]
$D(\vphi,\xi)=||\vphi-\xi||_0$,\ 
$\vphi,\xi\in E_0$;
\een
clearly, $(E_0,D)$ is a complete metric space.

Let $c\in I$  be fixed in the sequel. 
The {\it Razumikhin class}
of functions in $E_0$ attached to $c$, is defined as
\ben
\item[] (b01)\ \ 
$\calr_c=\{\vphi\in E_0; ||\vphi||_{0}=||\vphi(c)||\}$.
\een

Note that $\calr_c$ is nonvoid;
because any constant function belongs to it.
To substantiate this assertion, a lot of
preliminary facts are needed.

{\bf (A)}
For each $u\in E$, let $H[u]$ 
denote the constant 
function of $E_0$, defined as
\ben
\item[] (b02)\ \  
$H[u](t)=u$,\ $t\in I$.
\een
Note that, by this definition,
$$
||H[u]||_{0}=||u||,\
H[u](c)=u;
$$
whence $H[u]\in \calr_c$.
Denote, for simplicity
$\calk=\{H[u]; u\in E\}$;
this will be referred to as the
{\it constant class} of $E_0$.
The following properties of this class
are almost immediate;
so, we do not give details.

\bprop \label{p1}
Under the above conventions,
\ben
\item[] (cc1)\ \ 
$H[u+v]=H[u]+H[v]$,\ $H[\lb u]=\lb H[u]$,\
$\forall u,v\in E$, $\forall \lb\in R$;\\ 
hence, $\calk$ is a linear subspace of $E_0$
\item[] (cc2)\ \ 
$||u||=||H[u]||_0$, $\forall u\in E$
\item[] (cc3)\ \ 
the mapping $u\mapsto H[u]$ is an algebraic 
and topological isomorphism between 
$(E,||.||)$ and $(\calk,||.||_{0})$
\item[] (cc4)\ \ 
$\calk$ is $D$-complete (hence, $D$-closed) in $E_0$.
\een
\eprop

{\bf (B)}
Returning to the general case,
the following simple property holds.

\blemma \label{le1}
The Razumikhin class $\calr_c$ is
homogeneous, in the sense
\beq \label{201}
\lb\calr_c= \calr_c,\
\forall \lb\in R\sm \{0\};\ 
\mbox{whence},\  \calr_c=-\calr_c.
\eeq
\elemma

\bproof
Given [$\lb\in R$, $\vphi\in \calr_c$],
denote $\xi=\lb\vphi$. 
By definition,
$$
||\xi||_0=|\lb|\cdot||\vphi||_0,\
||\xi(c)||=|\lb|\cdot||\vphi(c)||.
$$
This, along with 
the choice of $\vphi$,
gives $\xi\in \calr_c$;
and completes the argument.
\eproof

{\bf (C)}
Let $\calt:E_0\to E$ be a 
(nonself) mapping.
We say that $\vphi\in E_0$ is a 
{\it PPF dependent fixed point} of $\calt$,
when $\calt\vphi=\vphi(c)$.
The class of all these will be denoted as
PPF-$\Fix(\calt;E_0)$.

Concerning existence and uniqueness 
properties involving such points,
the first contribution to this theory is
the 1977 statement in 
Bernfeld et al 
\cite{bernfeld-lakshmikantham-reddy-1977};\
referred to as: BLR theorem.
The following concepts and constructions are
necessary.

{\bf I)}
Given $k\ge 0$, call $\calt$, {\it $k$-contractive}, provided
\ben
\item[] (b03)\ \ 
$d(\calt\vphi,\calt\xi)\le k D(\vphi,\xi)$,\ 
for all $\vphi,\xi\in E_0$.
\een

{\bf II)}
Let us introduce a relation $(\nabla)$ 
over $E_0$, according to:
\ben
\item[] (b04)\ \ 
$\vphi\nabla \xi$ iff
$\calt\vphi=\xi(c)$ and $\vphi-\xi\in \calr_c$.
\een

{\bf III)}
Finally, for the starting element
$\vphi_0\in E_0$,
let us say that the sequence 
$(\vphi_n; n\ge 0)$ in $E_0$ is
{\it $(\vphi_0,\nabla)$-iterative}, in case
$\vphi_n \nabla \vphi_{n+1}$,\ $\forall n$.

\btheorem \label{t2}
Suppose that $\calt$ is $k$-contractive, for some
$k\in [0,1[$.
Then,

{\bf i)}
Given the starting point $\vphi_0\in E_0$, any 
$(\vphi_0,\nabla)$-iterative sequence
$(\vphi_n; n\ge 0)$ in $E_0$,
$D$-converges to an element of PPF-$\Fix(\calt;E_0)$.

{\bf ii)}
Given the couple of starting points 
$\vphi_0,\xi_0\in E_0$,
and letting $(\vphi_n; n\ge 0)$, $(\xi_n; n\ge 0)$
be a
$(\vphi_0,\nabla)$-iterative sequence
and
$(\xi_0,\nabla)$-iterative sequence, respectively,
we have the evaluation, for all $n\ge 0$,
\beq \label{202}
D(\vphi_n,\xi_n)\le (1/(1-k)) 
[D(\vphi_0,\vphi_1)+D(\xi_0,\xi_1)]+D(\vphi_0,\xi_0).
\eeq
In particular, if $\vphi_0=\xi_0$,
we have for all $n\ge 0$,
\beq \label{203}
D(\vphi_n,\xi_n)\le (2/(1-k)) D(\vphi_0,\vphi_1).
\eeq

{\bf iii)}
Let $(\vphi_n; n\ge 0)$, $(\xi_n; n\ge 0)$ be as above.
If $\vphi_n-\xi_n\in \calr_c$, for all $n\ge 0$,
then, necessarily, 
$\lim_n\vphi_n=\lim_n\xi_n$.

{\bf iv)}
If $\vphi^*, \xi^*$ 
are in PPF-$\Fix(\calt;E_0)$,
and $\vphi^*-\xi^*\in \calr_c$,
then $\vphi^*=\xi^*$.
\etheorem

[For completeness reasons, we shall provide 
a proof of the above result, at the end of this exposition].
\sk

{\bf (D)}
Technically speaking, the BLR theorem is conditional 
in nature; because, for the starting
$\vphi_0\in \calr_c$
[hence, all the more, 
for the starting 
$\vphi_0\in E_0$], 
the set of all $(\vphi_0,\nabla)$-iterative sequence
$(\vphi_n; n\ge 0)$ in $E_0$
(taken as before) may be empty. 
To avoid this drawback, we have two possibilities:

{\bf Option-1)}
all considerations above 
are to be restricted
to the constant class $\calk\incl \calr_c$;
which, as a $D$-closed linear 
subspace of $E_0$, yields an appropriate  
setting for any algebraic and/or topological reasoning 
to be applied

{\bf Option-2)}
the initial Razzumikhin class 
$\calr_c$ remains as it stands; 
but, with the price of imposing
further (strong) restrictions upon it.

A discussion of these is to be sketched 
under the lines below.

{\bf Part 1.}
Concerning the first option,
we note that, the imposed
$k$-contractive condition 
upon $\calt$ relates elements
$\vphi,\xi\in E_0$ with elements
$\calt\vphi,\calt\xi\in E$.
On the other hand, at the level 
of constant class $\calk$,
the underlying condition writes
\ben
\item[] (b05)\ \ 
$||\calt\vphi-\calt\xi||\le k ||\vphi(c)-\xi(c)||$; 
for all $\vphi,\xi\in \calk$;
\een 
so that, it relates elements $\vphi(c),\xi(c)\in E$  
with elements $\calt\vphi,\calt\xi\in E$.
Hence,
as long as we have a selfmap of $E$
that sends $\psi(c)\in E$ to $\calt\psi\in E$
(for all $\psi\in \calk$),
the last condition is of selfmap type.
The effectiveness of such a
construction is illustrated
by the considerations below.
Let $T:E\to E$ be the selfmap of $E$ introduced as
\ben
\item[] (b06)\ \ 
$Tu=\calt(H[u])$,\ \ $u\in E$.
\een

\bprop \label{p2}
Under these conventions, 
the following are valid:

{\bf i)}
If $x\in E$ is a fixed point of $T$, 
then
$\xi:=H[x]\in \calk$ is a 
PPF dependent fixed point of $\calt$

{\bf ii)}
Conversely, if $\zeta=H[z]\in \calk$ is a 
PPF dependent fixed point of $\calt$, then
$z\in E$ is a fixed point of $T$.
\eprop

\bproof
i)
If $x\in E$ is a fixed point of $T$, we have
$x=Tx=\calt(H[x])$;
or, equivalently,
$$
\mbox{
$x=\calt(\xi)$,\ where $\xi:=H[x]\in \calk$.
}
$$
This, by definition, gives
$\xi(c)=\calt(\xi)$;
which tells us that 
$\zeta\in \calk$ is a 
PPF dependent fixed point of $\calt$.

ii)
Suppose that 
$\zeta=H[z]\in \calk$ is a 
PPF dependent fixed point of $\calt$;
that is:
$\zeta(c)=\calt(\zeta)$.
This yields (by these notations)
$$
z=\calt(H[z])=Tz;
$$
so that, $z\in E$ is a fixed point of $T$.
\eproof

Having these precise, 
we may now proceed to the 
formulation of announced result.
Two basic concepts appear.

{\bf I)}
Given $k\ge 0$, call $\calt$, {\it $k$-contractive}, provided
\ben
\item[] (b07)\ \ 
$d(\calt\vphi,\calt\xi)\le k D(\vphi,\xi)$,\ 
for all $\vphi,\xi\in E_0$.
\een
In particular, with $\vphi=H[x]$, $\xi=H[y]$
(where $x,y\in E$), this relation becomes
\ben
\item[] (b08)\ \ 
$d(Tx,Ty)\le k D(H[x],H[y])=kd(x,y)$,\ 
for all $x,y\in E$;
\een
or, in other words:
the associated selfmap $T:E\to E$ is 
$(d;k)$-contractive.

{\bf II)}
For the arbitrary fixed $\vphi_0=H[x_0]\in \calk$,
we say that the sequence $(\vphi_n=H[x_n]; n\ge 0)$ 
in $\calk$ is {\it $(\vphi_0,\calt)$-iterative},
provided 
[$\calt\vphi_n=\vphi_{n+1}(c)$, $\forall n\ge 0$].
Note that, the family of such sequences is nonempty. 
In fact, for the starting $x_0\in E$,
the $(x_0,T)$-iterative sequence
$(x_n; n\ge 0)$ in $E$ is
well defined, according to 
the formula $(Tx_n=x_{n+1}; n\ge 0)$.
But then, under the above notations,
$$
\calt\vphi_n=Tx_n=x_{n+1}=\vphi_{n+1}(c),\
\forall n\ge 0;
$$
or, in other words: the sequence  
$(\vphi_n:=H[x_n]; n\ge 0)$ in $\calk$ 
is 
$(\vphi_0,\calt)$-iterative.
Conversely, 
if the sequence 
$(\vphi_n:=H[x_n]; n\ge 0)$ in $\calk$
is 
$(\vphi_0,\calt)$-iterative, then
(by the same formula),
the sequence $(x_n; n\ge 0)$ 
of $E$ is $(x_0,T)$-iterative;
hence, 
for each $\vphi_0\in \calk$,
the family of all
$(\vphi_0,\calt)$-iterative sequences 
in $\calk$ is a singleton. 
  
Putting these together, it 
follows, from
the Banach fixed point principle
we just exposed,
the following 
coincidence point result
involving our data
(referred to as:
Constant BLR theorem):

\btheorem \label{t3}
Suppose that 
$\calt$ is $k$-contractive, for some $k\in [0,1[$.
Then, 

{\bf i)}
$\calt$ has (in $\calk$) a unique 
PPF dependent fixed point  $\vphi^*$ in  $\calk$

{\bf ii)}
for the arbitrary fixed $\vphi_0\in \calk$,
the $(\vphi_0,\calt)$-iterative sequence
$(\vphi_n; n\ge 0)$ in $\calk$, 
$D$-converges  to $\vphi^*$.
\etheorem

Note that -- 
unlike the situation encountered at BLR theorem --
the iterative sequences above are 
constructible in a precise way 
[by means of the associated selfmap $T$];
so, this result is an effective one.

{\bf Part 2.}
The second option above 
starts from the fact that,
the construction in $\calr_c$ of iterative sequences
given by BLR theorem
requires the structural condition 
\ben
\item[] (b09)\ \ 
$\calr_c$ is  {\it algebraically closed}: 
$\vphi,\xi\in \calr_c \limpl \vphi-\xi \in\calr_c$
\een
(also referred to as:
$\calr_c$ is {\it a-closed});
this assertion seems to have been 
formulated, for the first time, in 
Dhage \cite[Observation I]{dhage-2012-JNSA}.
On the other hand, the existence 
property above is retainable, 
at the level of $\calr_c$,
when  
\ben
\item[] (b10)\ \ 
$\calr_c$ is topologically closed:
$\calr_c$ is a $D$-closed part of $E_0$;
\een
cf. 
Dhage \cite[Observation II]{dhage-2012-JNSA}.
Summing up, the following
"existential" version of Theorem \ref{t2} 
(referred to as: Existential BLR theorem)
enters into our discussion:

\btheorem \label{t4}
Suppose that 
$\calt$ is $k$-contractive, for some $k\in [0,1[$.
In addition, let us assume that $\calr_c$ is 
algebraically and topologically closed.
Then, $\calt$ has a unique 
PPF dependent fixed point in $\calr_c$.
\etheorem

Note that, 
this result is {\it not}
present in the 1977 paper by
Bernfeld et al 
\cite{bernfeld-lakshmikantham-reddy-1977}.
The above formulation is 
a quite recent  "by-product" of 
BLR theorem, 
under the lines imposed by the above remarks;
cf.
Kutbi and Sintunavarat
\cite{kutbi-sintunavarat-2014}.

Concerning the structural requirements above,
we stress that, from a methodological perspective,
the algebraically closed condition 
is a very strong one.
Before explaining our assertion,
let us give a useful characterization of this
concept.
[Since the verification is
almost immediate, we omit the details].

\blemma \label{le2}
The following conditions are equivalent:
\ben
\item[] (b11)\ \ 
$\calr_c$ is algebraically closed
\item[] (b12)\ \ 
$\calr_c$ is additive:
$\vphi,\psi\in \calr_c$ 
$\limpl$ $\vphi+\psi\in \calr_c$
\item[] (b13)\ \ 
$\calr_c$ is a linear subspace of $E_0$.
\een
\elemma

Having these precise, 
we are in position to
motivate our previous affirmation.

\bprop \label{p3}
Suppose that $\calr_c$ is algebraically closed;
or, equivalently, additive.
Then, 
necessarily, $\calr_c=\calk$;
hence $\calr_c$ is topologically closed as well.
\eprop

\bproof
Suppose that $\calr_c\sm \calk\ne \es$;
and take some function
$\vphi$ in this set difference; 
hence, in particular,
\beq \label{204}
\mbox{
$\vphi(r)\ne \vphi(c)$,
for at least one $r\in I$.
}
\eeq
The function $\xi:=H[\vphi(c)]$
belongs to the constant 
class $\calk$;
hence, to the Razumikhin class $\calr_c$
as well.
As $\calr_c$ is algebraically closed,
the difference function $\de:=\vphi-\xi$ is 
an element of $\calr_c$;
so that 
$$
||\de||_0=||\de(c)||=0;
$$
or, equivalently,
$$
\vphi(t)=\vphi(c),\ \forall t\in I;
$$
in contradiction with  
the initial choice of $\vphi$.
Hence, $\calr_c=\calk$, as claimed.
The last affirmation is immediate, by 
the topological properties of $\calk$
(see above). 
\eproof

Summing up,
the following conclusions are 
to be noted:
\sk

{\bf Conc-1)}
If
the Razumikhin class $\calr_c$ 
is algebraically closed,
we must have $\calr_c=\calk$;
so that,
Existential BLR theorem
over $\calr_c$ is identical --
in a trivial way --
with Constant BLR theorem
over $\calk$;
which -- 
as already shown --  
reduces to the 
Banach fixed point principle,
without imposing any regularity 
condition upon $\calr_c$.
This means that,
the algebraic (and/or topological)
closeness condition 
upon $\calr_c$ 
has a null generalizing effect, 
relative to 
constant BLR theorem.

{\bf Conc-2)}
If
the Razumikhin class $\calr_c$ 
is algebraically closed,
we have 
[in view of $\calr_c=\calk$], 
that any PPF fixed point result
on $\calr_c$ is identical --
in a trivial way --
with the corresponding 
PPF fixed point statement
over $\calk$;
which 
[by the associated selfmap construction] 
is reducible to 
(standard) fixed point theorems
over the (complete) metric space $(E,d)$,
without imposing any regularity 
condition upon $\calr_c$.
So, as before,
the algebraic (and/or topological) 
closeness condition upon $\calr_c$ 
has a null generalizing effect
relative to the 
(variant of considered result over)
constant class $\calk$.
\sk

In particular, this latter conclusion
tells us that   
all recent PPF dependent fixed point results --
based on such conditions upon $\calr_c$ --
obtained in     
Ciri\'{c} et al \cite{ciric-alsulami-salimi-vetro-2014}.
Dhage \cite{dhage-2012-FPT, dhage-2012-JNSA},
Harjani et al \cite{harjani-lopez-sadarangani-2010},
Hussain et al \cite{hussain-khaleghizadeh-salimi-akbar-2013},
Kaewcharoen \cite{kaewcharoen-2013},
Kutbi and Sintunavarat \cite{kutbi-sintunavarat-2014},
are in fact reducible 
to PPF dependent fixed point results over 
the constant class $\calk$;
and these, in turn, are deductible 
from
corresponding fixed point theorems 
over the supporting 
metrical structure $(E,d)$,
without imposing any regularity 
condition to $\calr_c$.
A verification of 
this assertion for the
PPF dependent fixed point result in
Agarwal et al \cite{agarwal-kumam-sintunavarat-2013}
is performed in the rest of our exposition.
The remaining cases are to be treated 
in a similar way; we do not give details.

\section{SVV type contractions}
\setcounter{equation}{0}

Let $X$ be a nonempty set.
Take a metric $d:X\times X\to R_+$ over it;
as well as a selfmap 
$T\in \calf(X)$.
The basic directions under which 
the question of determining the fixed points of $T$
is to be solved were already sketched.
As precise, a classical result in this direction is
the 1922 one due to 
Banach \cite{banach-1922}.
In the following, a conditional version 
of the quoted statement is given, 
under the lines proposed in
Samet et al \cite{samet-vetro-vetro-2012}.
Let the mapping $\al:X\times X\to R_+$
be fixed in the sequel.

{\bf I)}
We say that $T$ is 
{\it $\al$-admissible}, if
\ben
\item[] (c01)\ \ 
($\forall x,y\in X$):\
$\al(x,y)\ge 1$ 
implies $\al(Tx,Ty)\ge 1$.
\een

{\bf II)}
Given $k\ge 0$, we say that $T$ is 
{\it (SVV type) $(\al,k)$-contractive}, provided
\ben
\item[] (c02)\ \ 
$\al(x,Tx)\al(y,Ty)
d(Tx,Ty)\le kd(x,y)$,\
for all $x,y\in X$.
\een

{\bf III)}
Further, let us say that 
$\al$ is {\it $X$-closed}, provided
\ben
\item[] (c03)\ \ 
whenever the sequence $(x_n; n\ge 0)$ 
in $X$ and the element $x\in X$
fulfill 
[$\al(x_n,Tx_n)\ge 1$, $\forall n$], and $x_n\dtends x$, then
$\al(x,Tx)\ge 1$.
\een

{\bf IV)}
Finally, let us say that $T$ is 
{\it $X$-starting}, if 
\ben
\item[] (c04)\ \ 
there exists $x_0\in X$,
such that
$\al(x_0,Tx_0)\ge 1$. 
\een

The following conditional fixed point result
(referred to as: SVV theorem)
involving these data is then available:

\btheorem \label{t5}
Suppose that $\calt$ is 
(SVV type) $(\al,k)$-contractive,
for some $k\in [0,1[$.
In addition, suppose that 
$(X,d)$ is complete,
$T$ is $\al$-admissible,
$\al$ is $X$-closed,
and
$T$ is $X$-starting.
Then,

{\bf i)}
For the arbitrary fixed $x_0\in X$
with $\al(x_0,Tx_0)\ge 1$,
the sequence $(x_n; n\ge 0)$ in $X$ 
defined as
$(x_{n+1}=Tx_n; n\ge 0)$
converges to a fixed point 
$x^*\in X$ of $T$, with
$\al(x^*,Tx^*)\ge 1$

{\bf ii)}
$T$ has exactly one 
fixed point $x^*\in X$,
such that $\al(x^*,Tx^*)\ge 1$.
\etheorem

\bproof
There are two assertions 
to be clarified.

{\bf Step 1.}
Let us firstly verify that 
$T$ cannot have more than 
one fixed point $x^*\in X$, such that 
$\al(x^*,Tx^*)\ge 1$.
In fact, assume that,
$T$ would have another 
fixed point $y^*\in X$ 
such that $\al(y^*,Ty^*)\ge 1$.
From the (SVV type) contractive condition, we have
$$
\al(x^*,Tx^*)\al(y^*,Ty^*)d(Tx^*,Ty^*)\le kd(x^*,y^*).
$$
This, along with the imposed properties, yields
$d(Tx^*,Ty^*)\le k d(x^*,y^*)$;
or, equivalently (as $x^*,y^*$ are fixed points)
$d(x^*,y^*)\le kd(x^*,y^*)$;
wherefrom (as $0\le k< 1$),
$d(x^*,y^*)=0$; hence, $x^*=y^*$.

{\bf Step 2.}
Let us now establish the existence part.
Fix in the following some 
$x_0\in X$ with $\al(x_0,Tx_0)\ge 1$;
and let $(x_n:=T^nx_0; n\ge 0)$ stand for the
iterative sequence generated by it.
From this very choice,
$$
\al(x_0,x_1)=\al(x_0,Tx_0)\ge 1;
$$
whence, by the admissible property of $T$,
$$
\al(x_1,Tx_1)=\al(x_1,x_2)=\al(Tx_0,Tx_1)\ge 1;
$$
and so on. 
By a finite induction procedure,
one gets an evaluation like
\beq \label{301}
\al(x_n,Tx_n)=\al(x_n,x_{n+1})\ge 1,\ \forall n.
\eeq 
This tells us that the (SVV type) contractive condition
applies to  $(x_n,x_{n+1})$, for all $n\ge 0$.
An effective application of it gives
$$
\al(x_n,x_{n+1})\al(x_{n+1},x_{n+2})
d(x_{n+1},x_{n+2})\le kd(x_{n},x_{n+1}),\ 
\forall n;
$$
wherefrom
$$
d(x_{n+1},x_{n+2})\le kd(x_{n},x_{n+1}),\ 
\forall n.
$$ 
This, again by a finite induction procedure, gives
$$
d(x_n,x_{n+1})\le k^n d(x_0,x_1),\ \forall n;
$$
so that, as the series 
$\sum_n k^n$ converges, $(x_n; n\ge 0)$ is 
$d$-Cauchy.
As $(X,d)$ is complete,
$x_n\dtends x^*$ for some $x^*\in X$;
moreover, combining with (\ref{301})
(and the closed property of $\al$), 
one derives $\al(x^*,Tx^*)\ge 1$.
The (SVV type) contractive condition is
applicable to each pair $(x_n,x^*)$, 
where $n\ge 0$; and gives 
$$
\al(x_n,Tx_n)\al(x^*,Tx^*)
d(Tx_n,Tx^*)\le kd(x_n,x^*),\ \forall n;
$$
wherefrom
\beq \label{302}
d(x_{n+1},Tx^*)=d(Tx_n,Tx^*)\le kd(x_n,x^*),\ \forall n.
\eeq
The sequence $(y_n:=x_{n+1}; n\ge 0)$
is a subsequence of $(x_n; n\ge 0)$;\
so that, $y_n\dtends x^*$ as $n\to \oo$.
Passing to limit in (\ref{302}), gives
$d(x^*,Tx^*)= 0$; wherefrom
(as $d$ is sufficient), 
$x^*=Tx^*$.
The proof is thereby complete.
\eproof

Note that, the obtained result is not 
very general in the area;
but, it will suffice for our purposes. 
Various extensions of it 
may be found in 
Samet et al \cite{samet-vetro-vetro-2012}.

\section{AKS theorem}
\setcounter{equation}{0}

Under these preliminaries,
we may now pass to
the announced result concerning
PPF dependent fixed points.
\sk

Let $(E,||.||)$ be a Banach space;
and let $d$ be the induced by norm metric on $E$
($d(x,y)=||x-y||$,\ $x,y\in E$);
hence, $(E,d)$ is complete.
Further, let $I=[a,b]$ be a closed real interval,
and $E_0:=C(I,E)$ stand for the class 
of all continuous functions 
$\vphi:I\to E$, endowed with the supremum norm
($||\vphi||_{0}=\sup\{||\vphi(t)||; t\in I\}$,\ $\vphi\in E_0$).
As before, let $D$ stand for the 
induced by norm metric on $E_0$
($D(\vphi,\xi)=||\vphi-\xi||_0$,\ $\vphi,\xi\in E_0$);
clearly, $(E_0,D)$ is complete.

Let $c\in I$  be fixed in the sequel. 
The {\it Razumikhin class}
of functions in $E_0$ attached to $c$, is defined as
\ben
\item[] 
$\calr_c=\{\vphi\in E_0; ||\vphi||_{0}=||\vphi(c)||\}$.
\een
Note that $\calr_c$ is nonvoid;
because any constant function belongs to it.
Precisely,
for each $u\in E$, let $H[u]$ 
stand for the constant 
function of $E_0$, defined as:
$H[u](t)=u$,\ $t\in I$.
Note that, by this definition,
$||H[u]||_{0}=||u||$, $H[u](c)=u$;
whence $H[u]\in \calr_c$.
Denote, for simplicity,
$\calk=\{H[u]; u\in E\}$;
this will be referred to as the
{\it constant class} of $E_0$.

Now, let $\calt:E_0\to E$ be a mapping.
We say that $\vphi\in E_0$ is a 
{\it PPF dependent fixed point} of $\calt$,
when $\calt\vphi=\vphi(c)$.
As already noted, the natural 
way to determine such points
is offered by the constant class
$\calk$ of $E_0$.
Then, the points in question 
appear as fixed points of the 
selfmap $T:E\to E$, introduced as
$$
Tu=\calt(H[u]),\ \ u\in E.
$$
Precisely (see above)

{\bf fp-1)}
If $z\in E$ is a fixed point of $T$, 
then
$\zeta:=H[z]\in \calk$ is a 
PPF dependent fixed point of $\calt$

{\bf fp-2)}
Conversely, if $\vphi=H[u]\in \calk$ is a 
PPF dependent fixed point of $\calt$, then
$u\in E$ is a fixed point of $T$.
\sk

The fixed point result to be applied is
SVV theorem;
so, we have to clarify whether
the required conditions in terms of $T$ 
are obtainable 
via (nonself type) hypotheses in terms of $\calt$.
For an easy reference, 
we list the latter conditions.
Let in the following 
$\al:E\times E\to R_+$ be a mapping.

{\bf I)}
We say that $\calt$ is 
{\it $\al$-admissible}, if
\ben
\item[] (d01)\ \ 
($\forall \vphi, \xi\in E_0$):\
$\al(\vphi(c),\xi(c))\ge 1$ 
implies $\al(\calt\vphi,\calt\xi)\ge 1$.
\een
In particular, take 
$\vphi=H[x]$, $\xi=H[y]$, where $x,y\in E$.
Then, this condition yields
\ben
\item[]  
($\forall x, y\in E$):\
$\al(x,y)\ge 1$ 
implies $\al(Tx,Ty)\ge 1$;
\een
or, in other words:
the associated selfmap $T:E\to E$ is
$\al$-admissible.

{\bf II)}
Given $k\ge 0$, we say that $\calt$ is 
{\it $(\al,k)$-contractive}, provided
\ben
\item[] (d02)\ \ 
$\al(\vphi(c),\calt\vphi)\al(\xi(c),\calt\xi)
d(\calt\vphi,\calt\xi)\le kD(\vphi,\xi)$,\
$\forall \vphi,\xi\in E_0$.
\een
As before, take 
$\vphi=H[x]$, $\xi=H[y]$, where $x,y\in E$.
Then, by this condition,
\ben
\item[] 
$\al(x,Tx)\al(y,Ty)d(Tx,Ty)\le kd(x,y)$,\
for all $x,y\in E$;
\een
i.e.: $T$ is (SVV type) $(\al,k)$-contractive.

{\bf III)}
Further, let us say that 
$\al$ is {\it $E_0$-closed}, provided
\ben
\item[] (d03)\ \ 
whenever the sequence $(\vphi_n; n\ge 0)$ 
in $E_0$ and the element $\vphi\in E_0$
fulfill 
[$\al(\vphi_n(c),\calt\vphi_n)\ge 1$, $\forall n$], 
and $\vphi_n\Dtends \vphi$, then
$\al(\vphi(c),\calt\vphi)\ge 1$.
\een
In particular, taking
$(\vphi_n=H[x_n]; n\ge 0)$ and 
$\vphi=H[x]$, for some 
sequence $(x_n; n\ge 0)$ in $E$ and some
point $x\in E$, this requirement becomes:
\ben
\item[]  
whenever the sequence $(x_n; n\ge 0)$ 
in $E$ and the element $x\in E$
fulfill [$\al(x_n,Tx_n)\ge 1$, $\forall n$], 
and $x_n\dtends x$, then
$\al(x,Tx)\ge 1$;
\een 
or, in other words: $\al$ is $E$-closed.

{\bf IV)}
Finally, let us say that $\calt$ is 
{\it $E_0$-starting}, if 
\ben
\item[] (d04)\ \ 
there exists $\vphi_0\in E_0$,
such that
$\al(\vphi_0(c),\calt\vphi_0)\ge 1$.
\een
Likewise, 
let us say that $\calt$ is 
{\it $\calk$-starting}, if 
\ben
\item[] (d05)\ \ 
there exists $\vphi_0\in \calk$,
such that
$\al(\vphi_0(c),\calt\vphi_0)\ge 1$. 
\een

Clearly, if
$\calt$ is 
$\calk$-starting, 
then $\calt$ is 
$E_0$-starting as well.
The reciprocal inclusion holds too,
under an admissible property upon $T$.

\bprop \label{p4}
Assume that 
the (nonself mapping) $\calt$ is $\al$-admissible 
and
$E_0$-starting.
Then, $\calt$ is $\calk$-starting.
\eprop

\bproof
As $\calt$ is $E_0$-starting,
there exists $\vphi_0\in E_0$
such that $\al(\vphi_0(c),\calt\vphi_0)\ge 1$.
On the other hand, for
the point $\calt\vphi_0\in E$, 
we may consider the 
function
$\xi_0=H[\calt\vphi_0]$ 
in the constant class 
$\calk$. 
This, by definition, means
$$
\xi_0(t)=\calt\vphi_0,\ \forall t\in I;\
\mbox{whence},\ 
\xi_0(c)=\calt\vphi_0. 
$$
The starting property of $\calt$ becomes:
$$
\al(\vphi_0(c),\xi_0(c))\ge 1.
$$
As $\calt$ is $\al$-admissible, 
this yields
$\al(T\vphi_0,T\xi_0)\ge 1$.
Combining with a preceding relation, 
we thus have
$$
\al(\xi_0(c),\calt\xi_0)\ge 1;
$$
which tells us that 
$\calt$ is $\calk$-starting.
\eproof

Now, as 
$\calt$ is $\calk$-starting,
there exists $\vphi_0=H[x_0]$
(where $x_0\in E$)
in $\calk$, such that
$\al(\vphi_0(c),\calt\vphi_0)\ge 1$. 
This, by definition,
means $\al(x_0,Tx_0)\ge 1$;
and tells us that 
the associated selfmap $T$ is 
$E$-starting. 
\sk

Putting these together,
it results, via 
SVV theorem
(and our preliminary facts),
the following 
PPF dependent fixed point result
(referred to as: AKS theorem),
involving these data:

\btheorem \label{t6}
Suppose that $\calt$ is
$(\al,k)$-contractive,
for some $k\in [0,1[$.
In addition, suppose that 
$\calt$ is $\al$-admissible,
$\al$ is $E_0$-closed, 
and
$\calt$ is $E_0$-starting.
Then

{\bf i)}
$\calt$ is $\calk$-starting (see above)

{\bf ii)}
For the arbitrary fixed $\vphi_0\in \calk$
with $\al(\vphi_0(c),T\vphi_0)\ge 1$,
the sequence $(\vphi_n; n\ge 0)$ in $\calk$ 
defined as
$(\vphi_{n+1}(c)=\calt\vphi_n; n\ge 0)$
converges to a PPF dependent fixed point 
$\vphi^*\in \calk$ of $T$, with 
$\al(\vphi^*(c),\calt\vphi^*)\ge 1$

{\bf iii)}
$\calt$ has a exactly one PPF dependent
fixed point $\vphi^*$ in the constant class $\calk$,
such that $\al(\vphi^*(c),\calt\vphi^*)\ge 1$.
\etheorem

In particular, when 
the Razumikhin class
$\calr_c$ is 
algebraically (as well as topologically) closed,
and the starting condition is taken as
\ben
\item[] (d06)\ \ 
$\calt$ is 
{\it $\calr_c$-starting}:\ 
there exists $\vphi_0\in \calr_c$,
with
$\al(\vphi_0(c),\calt\vphi_0)\ge 1$, 
\een
the obtained result
is just the PPF dependent fixed point result in
Agarwal et al \cite{agarwal-kumam-sintunavarat-2013}.
However, as noted earlier, 
all these conditions have no
generalizing effect
upon our data;
because, in view of $\calr_c=\calk$,  
the obtained product is again AKS theorem;
which is deductible without any 
regularity condition upon $\calr_c$.

Finally, note that, by the same techniques,
it follows that
all PPF dependent fixed point results in 
Ciri\'{c} et al \cite{ciric-alsulami-salimi-vetro-2014},
Dhage \cite{dhage-2012-FPT, dhage-2012-JNSA},  
Harjani et al \cite{harjani-lopez-sadarangani-2010},
Hussain et al \cite{hussain-khaleghizadeh-salimi-akbar-2013},
Kaewcharoen \cite{kaewcharoen-2013}, 
Kutbi and Sintunavarat \cite{kutbi-sintunavarat-2014},
are in fact reducible to (corresponding)
fixed point statements over standard metric structures.
Some other aspects will be delineated elsewhere.

\section{Proof of BLR theorem}
\setcounter{equation}{0}

Let us now return to 
the BLR theorem.
For completeness reasons,
we shall provide a proof
of this result which, in part,
differs from the original one.

Let $(E,||.||)$ be a Banach space;
and let $d$ be the induced by norm metric on $E$
($d(x,y)=||x-y||$,\ $x,y\in E$);
hence, $(E,d)$ is complete.
Further, let $I=[a,b]$ be a closed real interval,
and $E_0:=C(I,E)$ stand for the class 
of all continuous functions 
$\vphi:I\to E$, endowed with the supremum norm
($||\vphi||_{0}=\sup\{||\vphi(t)||; t\in I\}$,\ $\vphi\in E_0$).
As before, let $D$ stand for the induced metric
($D(\vphi,\xi)=||\vphi-\xi||_0$,\ $\vphi,\xi\in E_0$);
clearly, $(E_0,D)$ is complete.

Let $c\in I$  be fixed in the sequel. 
The {\it Razumikhin class}
of functions in $E_0$ attached to $c$, is defined as
\ben
\item[] 
$\calr_c=\{\vphi\in E_0; ||\vphi||_{0}=||\vphi(c)||\}$.
\een
Note that $\calr_c$ is nonvoid;
because any constant function belongs to it.

Finally, let $\calt:E_0\to E$ be a (nonself) mapping.
We say that $\vphi\in E_0$ is a 
{\it PPF dependent fixed point} of $\calt$,
when $\calt\vphi=\vphi(c)$.
The class of all these will be denoted as
PPF-$\Fix(\calt;E_0)$.

To establish the existence and uniqueness
result in question,  
the following concepts and constructions are
necessary.

{\bf I)}
Given $k\ge 0$, call $\calt$, {\it $k$-contractive}, provided
\ben
\item[]  
$d(\calt\vphi,\calt\xi)\le k D(\vphi,\xi)$,\ 
for all $\vphi,\xi\in E_0$.
\een
Note that, as a consequence of this, $\calt$ is continuous:
$\vphi_n\Dtends \vphi$ implies 
$\calt\vphi_n\dtends \calt\vphi$.

{\bf II)}
Let us introduce a relation $(\nabla)$ 
over $E_0$, according to:
\ben
\item[]
$\vphi\nabla \xi$ iff
$\calt\vphi=\xi(c)$ and $\vphi-\xi\in \calr_c$.
\een

{\bf III)}
Finally, for the starting element
$\vphi_0\in E_0$,
let us say that the sequence 
$(\vphi_n; n\ge 0)$ in $E_0$ is
{\it $(\vphi_0,\nabla)$-iterative}, in case
$\vphi_n \nabla \vphi_{n+1}$,\ $\forall n$.
 
Having these precise, we may now proceed 
to the announced

\bproof {\bf (BLR theorem)}
There are several steps to be considered.

{\bf Part 1.}
Take a starting point $\vphi_0\in E_0$, 
and let
$(\vphi_n; n\ge 0)$ in $E_0$
be some
$(\vphi_0,\nabla)$-iterative sequence.
By definition, 
$\vphi_0\nabla \vphi_1$;
which means
$$
\calt\vphi_0=\vphi_1(c),\ 
d(\vphi_0(c),\vphi_1(c))=D(\vphi_0,\vphi_1).
$$
Further, $\vphi_1\nabla \vphi_2$;
which means
$$
\calt\vphi_1=\vphi_2(c),\ 
d(\vphi_1(c),\vphi_2(c))=D(\vphi_1,\vphi_2).
$$
Note that, as a combination of these, we have
(by the contractive property of $T$)
$$
D(\vphi_1,\vphi_2)= 
d(\vphi_1(c),\vphi_2(c))=
d(\calt\vphi_0,\calt\vphi_1)\le k D(\vphi_0,\vphi_1).
$$
This procedure may continue indefinitely;
and gives the iterative relations
\beq \label{501}
D(\vphi_{n+1},\vphi_{n+2})\le 
k D(\vphi_n,\vphi_{n+1}),\ \forall n\ge 0.
\eeq
By a finite induction procedure, one gets
\beq \label{502}
D(\vphi_{n},\vphi_{n+1})\le 
k^n D(\vphi_0,\vphi_1),\ \forall n\ge 0;
\eeq
and since the series $\sum_n k^n$ converges,
the sequence
$(\vphi_n; n\ge 0)$ is $D$-convergent.
As $(E_0,D)$ is complete,
$\vphi_n\Dtends \vphi^*$ as $n\to \oo$, for 
some $\vphi^*\in E_0$.
Passing to limit as $n\to \oo$
in the iterative relation that defines the
sequence $(\vphi_n; n\ge 0)$, one gets
(as $\calt$ is $(D,d)$-continuous) 
$\calt\vphi^*=\vphi^*(c)$;
i.e.; $\vphi^*$ is an element of
PPF-$\Fix(\calt;E_0)$.

{\bf Part 2.}
Take a couple of starting points 
$\vphi_0,\xi_0\in E_0$,
and let $(\vphi_n; n\ge 0)$, $(\xi_n; n\ge 0)$
be their corresponding
$(\vphi_0,\nabla)$-iterative sequence
and
$(\xi_0,\nabla)$-iterative sequence, respectively.
By the preceding part, 
$(\vphi_n; n\ge 0)$ fulfills the 
iterative relations (\ref{501}) and (\ref{502}).
Likewise, $(\xi_n; n\ge 0)$ fulfills
the iterative relations
\beq \label{503}
D(\xi_{n+1},\xi_{n+2})\le 
k D(\xi_n,\xi_{n+1}),\ \forall n\ge 0;
\eeq
wherefrom (by a finite induction procedure)
\beq \label{504}
D(\xi_{n},\xi_{n+1})\le 
k^n D(\xi_0,\xi_1),\ \forall n\ge 0.
\eeq
From the triangle inequality, 
one gets, for all $n\ge 1$,
$$
\barr{l}
D(\vphi_n,\xi_n)\le 
D(\vphi_{n-1},\vphi_n)+D(\xi_{n-1},\xi_n)+
D(\vphi_{n-1},\xi_{n-1})\le \\
k^{n-1}[D(\vphi_0,\vphi_1)+D(\xi_0,\xi_1)]+
D(\vphi_{n-1},\xi_{n-1}).
\earr
$$ 
By repeating the procedure, it follows that,
after $(n-1)$ steps, 
the first part of conclusion {\bf ii)} follows.
Moreover, if $\vphi_0=\xi_0$, then
$\calt\vphi_0=\calt\xi_0$
(that is: $\vphi_1(c)=\xi_1(c)$);
wherefrom (by the choice of our
iterative sequences)
$$
\barr{l}
D(\vphi_0,\vphi_1)=
d(\vphi_0(c),\vphi_1(c)), \\
D(\xi_0,\xi_1)=
d(\xi_0(c),\xi_1(c))=
d(\vphi_0(c),\vphi_1(c)); \\
\earr
$$
and, this yields 
the second part of conclusion {\bf ii)}.

{\bf Part 3.}
Let $(\vphi_n; n\ge 0)$, $(\xi_n; n\ge 0)$ be taken 
as in the premise above.
By the very definition of these points
and the working condition, we have, 
for all $n\ge 1$,
$$
D(\vphi_n,\xi_n)= 
d(\vphi_n(c),\xi_n(c))=
d(\calt\vphi_{n-1},\calt\xi_{n-1})\le \\
k D(\vphi_{n-1},\xi_{n-1}).
$$
This, by a finite induction procedure, gives 
$$
D(\vphi_n,\xi_n)\le 
k^n D(\vphi_0,\xi_0),\ \forall n.
$$
Passing to limit as $n\to \oo$,
and noting that both $(\vphi_n; n\ge 0)$ and
$(\xi_n; n\ge 0)$ are $D$-convergent 
(see above), conclusion {\bf iii)} follows.

{\bf Part 4.}
Let $\vphi^*, \xi^*$ 
be two elements in PPF-$\Fix(\calt;E_0)$,
taken as in the stated premise.
By the working condition and the 
contractive property of $T$,
$$
D(\vphi^*,\xi^*)=d(\vphi^*(c),\xi^*(c))
=d(\calt\vphi^*,\calt\xi^*)\le k D(\vphi^*,\xi^*);
$$
and, from this (as  $D$ is a metric), 
$\vphi^*=\xi^*$.
The proof is thereby complete.
\eproof

As above said, this result --
as well as its extensions due to
Pathak \cite{pathak-1993}
and
Som \cite{som-1985} --
is just a conditional one;
because the iterative sequences appearing there 
are not effective.
The correction of this drawback
by taking $\calr_c$ as algebraically closed
is, ultimately, without effect;
for, as the preceding developments  
show, it brings the PPF dependent fixed point problem 
in question at the level of constant class, $\calk$.
From a theoretical viewpoint,
this reduced problem may be of some avail 
in solving the (metrical) question we deal with.
However, from a practical perspective,
the resulting constant solutions for the 
nonlinear functional differential/integral 
equations based on these techniques --
such as, the ones in 
Kutbi and Sintunavarat
\cite{kutbi-sintunavarat-2014} --
are not very promising.
Further aspects will be discussed elsewhere.


\end{document}